\documentclass{article}
\usepackage{mflogo}
\usepackage{amsmath}
\usepackage{amsfonts}
\usepackage{hyperref}
\usepackage{amssymb}
\usepackage{graphics}
\newtheorem{theorem}{Theorem}

\title{Boundary behaviour of Loewner Chains}
\author{Alexander Kuznetsov}

\begin{document}
\maketitle
\begin{abstract}
 In paper found conditions that guarantee that solution of Loewner-Kufarev 
equation maps unit disc onto domain with quasiconformal rectifiable boundary, or it has continuation on closed unit disc, or 
it's inverse function has continuation on closure of domain.
\end{abstract}
\section{Introduction}
In recent time there was a significant interest to evolution processes of domains in complex 
plane. This phenomenas are well described by  equations of Loewner-Kufarev type. 

Let $D(t)$ be a simple connected domain in the complex plane, changing in time. Without loss of 
generality we can assume that $D(t)$ contains origin. Let a function $f(z,t)$ 
($f'(0,t)>0,f(0,t)=0$) for 
each fixed $t$ map unit disc  $\mathbb D =\{z:|z|<1\}$ onto domain $D(t)$. If  
$D(t_1)\subset 
D(t_2),t_1 > t_2$, then function $f(z,t)$ satisfies the Loewner-Kufarev equation

\begin{equation}
\label{partialLoewnerKufarev}
\frac{\partial f}{\partial t}=-\frac{\partial f}{\partial t}z p(z,t),z\in\mathbb{D},
\end{equation}
where $p(z,t):\mathbb D\times[0,+\infty)\to \mathbb C $ is such that for each fixed  $t\in[0,+\infty)$, $\Re e 
p(z,t) > 0,p(0,t)>0$ and for each fixed  $z\in\mathbb D$,  $p(z,t)$ is mesarable on interval
$[0,+\infty)$. Making change of variable~(\ref{partialLoewnerKufarev}) we can assume, that for 
each fixed $t$ function $p(z,t)$  belongs to class $\cal P$, consisting of functions
$p(0)=1,\Re e p(z)>0$.

Such type equations arise in description of Hele-Show 
flow~\cite{Galin,Polubarinova-Kochina1,Polubarinova-Kochina2,VasilevHS}. In 
paper~\cite{Carleson} L. Carleson and N. Makarov used this equations for the study of DLA  processes. 
Parametric method of representations of univalent functions~\cite{Aleksandrov,Pommerenke} is based on this approach. The 
importance of this method and approach is shown by L. de Branges proof of famous Bieberbach 
hypothesis~\cite{Branges}, that was the central problem of functions theory for long time.

If $p(z,t)=\frac{e^{iu(t)}+z} {e^{iu(t)}-z}$, where $u(t)$ is real function, then equation will be Loewner equation
~\cite{Loewner}. Recently there was a great activity in the study of this equation~\cite{Schramm1,Schramm2,Shramm3,Shramm4}, 
there $u(t)$ describes one-dimensional Brownian motion. In this case equation~(\ref{partialLoewnerKufarev}) 
becomes Stochastic Loewner Equation  (SLE). The study of geometrical characteristics of solutions of SLE,  by various 
parameters of Brownian motion is of a great interest.

In paper~\cite{Rohde} it was proved, that solution of Loewner equation maps unit disc on quasislit-disc if and only 
if function  $u(t)$ is Hoelder continuous with exponent 1/2.
Also in it was shown if equation 
$$\inf_{\varepsilon>0}\sup\limits_{|t-s|<\varepsilon}\frac{|u(t)-u(s)|}{\sqrt{|t-s|}}<c,$$
holds for some small $c>0$,( in paper of Lind~\cite{Lind} was proved, that $c=4$), when solution of Loewner equation
maps unit-disc on quasislit-disc.

Let ${\cal P}_{s}$ be some some subclass of class $\cal P$ and ${\cal PT}_{s}$ be a class of functions
$p(z,t):\mathbb D\times[0,+\infty)\to \mathbb C$ such that for each fixed  $t\in[0,+\infty)$, 
$p(z,t)\in{\cal P}_{s}$ and, for each fixed $z\in\mathbb D$,  $p(z,t)$ is mesarable on interval
$[0,+\infty)$. The main  goal of this paper is to find conditions of class ${\cal P}_{s}$, under that the
solution of equation~(\ref{partialLoewnerKufarev}) with $p(z,t)\in{\cal PT}_{s}$ will have certain boundary behaviour.

Similar question was studied by J. Becker. In papers~\cite{Becker1,Becker2} it was proved, if 

$$\Big|\frac{p(z,t)-1}{p(z,t)+1}\Big|\leqslant k<1,$$ 
then solution of equation~(\ref{partialLoewnerKufarev}) has 
quasiconformal extension. Recently A. Vasil'ev~\cite{Vasilev} gave a description of $p(z,t)$ under conditions that,  
solution $f(z,t)$ of~(\ref{partialLoewnerKufarev}) has quasiconformal extension for each $t$.

In the papers~\cite{Goryainov,Gutlyanski} were discussed problems of a similar nature, but from the point of view of 
covering and distortion theorems for various classes of univalent functions, corresponding to
various classes ${\cal P}_s$.

\section{Behaviour of inverse function}

Let us show that the study of boundary behaviour of solution of equation~(\ref{partialLoewnerKufarev}) for various classes
${\cal PT}_{s}$, can be reduced to the study of Loewner-Kufarev ODE.

Let us consider a function $w^*(z,t,s)=f^{-1}(f(z,t_2),t_1)$, that satisfies the following equation

$$
\frac{d w^*}{d t_1}=w^*p(w^*,t_1),\ w(z,s,s)=z,\, z\in\mathbb{D}.
$$

Let $t_1$ change from  $0$ to $t$, and $t_2=t$. When $w(z,0,t)=f(z,t)$. Assuming that $\tau=1-t_1$, we have,
that $w(z,t)=f(z,t)$, where $w(z,\tau)$ is solution of equation
\begin{equation}
\label{LoewnerKufarev}
\frac{dw}{d\tau}=-wp(w,\tau),w(z,0)=z,\,z\in\mathbb D.
\end{equation}
We have the necessary result.

\begin{theorem}
\label{inversContinuity} If there are constants $0< \alpha < 1$ and $C_1,C_2 > 0$ such that

\begin{equation}
\label{est1}
\frac{|p(z,t)|}{\Re e p(z,t)}\leqslant\frac{C_1}{(1-r)^\alpha},\ r=|z|,
\end{equation}

\begin{equation}
\label{est2}
\Re e p(z,t) \geqslant C_2(1-r)^\alpha,\ r=|z|.
\end{equation}
Then inverse function to solution $w(z,t)$ of equation~(\ref{LoewnerKufarev}) has a continuation on boundary of domain $w(\mathbb 
D,t)$.
\end{theorem}
Let us consider a integral curve $s_{(z_0,t_0)}$, passing through point $(z_0,t_0)$, and find estimate on
it's part, laying in cylinder $C_r=\{z:1/2<r<|z|<1\}$. 

Separating the real and the imaginary part in equation~(\ref{LoewnerKufarev}), we have

\begin{equation}
\label{separationRealAndImage}
	\begin{cases}
		\frac{d\rho}{dt}=-\rho\,\Re e p(w,t),\\
		\frac{d\arg w}{dt}=-\Im m p(w,t),
	\end{cases}
\end{equation}
where $\rho=|w|$.

From system~(\ref{separationRealAndImage}) follows that function  $\rho(t) = |w(z,t)|$ is monotonic 
on $t$. Thus, we can make a change of variable in equation~(\ref{LoewnerKufarev}) from $t$ on $\rho$. Making it, 
we have that length of part of integral curve $s_{(z_0,t_0)}$, that  laying in cylinder  $C_r$, equals to

$$
l_s(r) = \int\limits_{r}^{1}
	\sqrt{
		\Big(\frac{dw}{d\rho}\Big)^2 + 
		\Big(\frac{dt}{d\rho}\Big)^2
	}
d\rho.$$ 
Using inequalities~(\ref{est1}),~(\ref{est2}) and $|w(z,t)|\leqslant z$ we obtain 

$$
l_s(r)\leqslant \int\limits_{r}^{1}
\sqrt{
	\frac{C_1^2} {(1-\rho)^{2\alpha}} + 
	\frac{1} {C_2^2\rho^2(1-\rho)^{2\alpha}}
}
d\rho
	\leqslant 
C_3(1-r)^{1-\alpha}.
$$ 
Where $C_3$ is some constant, that depends only on  $C_1,C_2$. 
Thus,
\begin{equation}
\label{curveLenght}
	\lim\limits_{r\to 1}l_s(r)=0.
\end{equation}

Let $w(z,t)$ be a solution of equation~(\ref{LoewnerKufarev}),  passing through  $(z_0,t_0)$ and 
$\tau$ be minimal number, for that $w(z,t)$ exits on interval $(\tau,t_0)$. Using monotony of $|w(z,t)|$ 
on $t$ and~(\ref{curveLenght}), we have that there is $\lim\limits_{t\to\tau}w(z,t)=\gamma(z_0,t_0)$ and
$|\gamma(z_0,t_0)|=1$

Let define a function $s_{t_1}(z_0,t_0),|z|\leqslant 1,0\leqslant t_1 \leqslant t_0$ in a following way. 
If $w(z,t)$ solution of equation (\ref{LoewnerKufarev}) with initial condition $(z_0,t_0)$ 
exists for $t_1$,  we assume that $s_{t_1}(z_0,t_0)=w(z,t_1)$. Let  $Q_{0}$ be a set of pairs $(z,t_0)$ 
for that this is true. Let, $s_{t_1}(z_0,t_0)=\gamma(z_0,t_0)$, if $z\in 
Q_1=\{(z_0,t_0):|z_0|<1,(z_0,t_0)\notin Q_0\}$ and if $z\in Q_2=\{(z_0,t_0):|z_0|=1\}$, then
$s_{t_1}(z_0,t_0) = (z_0,t_0)$. 

Let us show that function $s_{t_1}(z_0,t_0)$ is a continuous one. 

By theorem of continuous depending of solution form initial conditions for differential equations  
in integral from (see for e.g.~\cite{Pontryagin} p.182),  
set $Q_0$ is open and function $s_{t_1}(z_0,t_0)$ is continuous on it. 

Let $(z_0,t_0)\in Q_2$ and point $(z',t')$ satisfy the following equation
$|z_0-z'|+|t_0-t'|<\varepsilon$. The length of part of curve $w(z,t)$, passing through point 
$(z',t')$ and laying in cylinder $C_\varepsilon$,  is not greater then  $C_3\varepsilon^{\alpha}$. So we have following 
inequality

$$|s_{t_1}(z_0,t_0)-s_{t_1}(z',t')|\leqslant C_3\varepsilon^{1-\alpha}+\varepsilon.$$ 
This inequality means that function $s_{t_1}(z,t)$ is continuous on set $Q_2$.

Let $(z_0,t_0)\in Q_1$,  $t^*$, such that $1-|w(z,t^*)|=2\varepsilon$ (for enough small  $\varepsilon$ 
such point exists) and $O$ be  neighbourhood of point  
$(t^*,w(z_0,t^*))$, that  is defined by the following inequality

\begin{equation}
\label{nebhInCaseOfQ1}
|t^*-t'|+|w(z,t^*)-z'|<\varepsilon. 
\end{equation}

Using the theorem of continuous depending of solution form initial conditions for differential equations we have, that 
there is  $\delta>0$ so that, form inequality $|t_0-t'_0|+|z_0-z'_0|<\delta$ follows that $w^*(z',t')\in 
O$, where  $w^*(z,t)$ integral curve passing through point  $(z'_0,t'_0)$. Taking into account, that the length of parts of 
integral curves $w^*(z,t)$ and $w(z,t)$, laying in cylinder $C_{2\varepsilon}$, is not greater then 
$2C_3(2\varepsilon)^{1-\alpha}$, we have

$$|s_{t_1}(z_0,t_0)-s_{t_1}(z',t')| < 2C_3(2\varepsilon)^{1-\alpha}+\varepsilon.$$ 
This inequality means that function $s_{t_1}(z,t)$ is continuous on set $Q_1$.
Noting that  $f^{-1}(z,t)=s_0(z,t)$ we have the necessary result.

\section{Boundary behaviour}
\begin{theorem}
\label{HolderTheorem} Following statmets are true:
\begin{enumerate}
 \item
If inequality
\begin{equation}
\label{inequalityForHoelder}
\frac{\Re e zp_z(z,t)}{|z|\Re e p(z,t)}\leqslant\frac{k}{1-r}+O\Big(\frac1{(1-r)^\alpha}\Big),\, k<1,r=|z| 
\end{equation}
holds, then function  $w(z,t)$ is Hoelder continuous with exponent  $1-k$ in closed unit disc.
 \item
If inequality
\begin{equation}
\label{inequalityForQuasiconforaml}
\frac{\Re e zp_z(z,t)}{|w|\Re e p(z,t)}=O\Big(\frac1{(1-r)^\alpha}\Big),\,r=|z| 
\end{equation} holds,
then function $w(z,t)$ maps unit disc on domain with rectifiable and quasiconformal boundary.
\end{enumerate} 

\end{theorem}

In the proof we are using ideas similar to ideas in papers~\cite{Goryainov,Gutlyanski}. 

Differentiating equation~(\ref{LoewnerKufarev}) on $z$, we obtain

$$
	\frac{dw_z}{dt} 
	= -w_z p(w,z) - w w_z p_z(w,t),
   |z| < 1 , w_z(z,0) = 1.
$$
Making change of variable from  $t$ to $\rho$, we have

$$
\frac{d w_z}
	  {d\rho} 
= - w_z 
	\frac{p(w,t)+w p_z(w,t)}
	{\rho\Re e p(w,t)}.
$$
Integrating of this equations gives us
$$
w_z(z,t) = e^{
		\int_{r_1}^{|z|}
			\frac{p + w p_z}{\rho\Re e p}
		d\rho
	}.
$$
From that we have 

\begin{equation}
\label{derivativeItegral}
|w_z(z,t)|
\leqslant e^{
		\int_{r_1}^{|z|}
			\frac1\rho+\frac{\Re e w p_z}{\rho\Re e p}
		d\rho
	}.
\end{equation}
Taking into account what $r_1=O(t)$ and inequality~(\ref{inequalityForHoelder}), we have

$$
C_1(1 - |z|)^k
	\leqslant
	|w_z(z,t)|
	\leqslant
C_2\frac{1}{(1 - |z|)^k}.
$$
From this follows that, function $w(z,t)$ is Hoelder continuous with exponent  $1-k$ in closed unit disc (see for 
e.g.~\cite{Pommerenke} p.300). 
First part of the theorem is proved.

Taking into account (\ref{inequalityForQuasiconforaml}) and (\ref{derivativeItegral}),  we have that, there is $t_0$ 
such that for all $t<t_0$ we have
$|w_z(z,t)-1|\leqslant 1/4$. Let us consider three points on circle $z_1,z_2,z_3$, 
so that point  $z_2$ lays between points $z_1,z_3$. 
From equality  $$|w(z_1,t)-w(z_2,t)|=\Big|\int_{z_1}^{z_2}w_z(z,t)dz\Big|,$$ we obtain

$$\frac34|z_1-z_2|\leqslant|w(z_1,t)-w(z_2,t)|\leqslant\frac54|z_1-z_2|.$$
So we have

$$\frac{|w(z_1,t)-w(z_2,t)|}{|w(z_1,t)-w(z_3,t)|}\leqslant\frac{5}{3}\frac{|z_1-z_2|}{|z_1-z_3|}.$$
Therefore, boundary of domain $w(\mathbb D,t)$ is quasiconformal curve and function  
$w(z,t)$ have quasiconformal extension on hole complex plane, for $t<t_0$.

Let $f_i(z)=w_i(z,t^*)$, where $t^*=t/n,i=1,2,\ldots,n$ and $w_i$ is solution of equation
$$
\frac{dw_i}{dt} = w_i p(w_i,t-i t^*), w_i(z,0)=z.
$$
When 

$$w(z,t)=f_1\circ f_2\circ\ldots \circ f_n.$$   

Thus choosing $n$ so that $t^*=t/n<t_0$ we have that all $f_i,i=1,\ldots,n$ have quasiconformal extension.
Therefore  $w(z,t)$ as well.
Also we have $|w_z(z,t)|\leqslant (3/2)^n$ so $w(z,t)$ maps unit disc on domain with rectifiable boundary. 
The second part of theorem is proved.

\section{Examples}
In this section  we consider various subclasses of $\cal P$, that have a simple geometrical description and satisfy conditions 
of theorems~\ref{inversContinuity} and \ref{HolderTheorem}.
Let $\Omega$ be convex subdomain of the right half-plan and 
${\cal P}_{\Omega}$ be class of functions $p:\mathbb D \to \Omega$. Through 
$p_\Omega$ we define a map from $\mathbb D$ onto $\Omega$. Any function  $p\in{\cal P}_{\Omega}$ can be represented as
$p_\Omega(\varphi(z))$, where $\varphi:\mathbb D \to \mathbb D$ and $\varphi(0)=0$
From this and  Schwarz lemma follows that if function $p_\Omega$ satisfies conditions of  
theorem~\ref{inversContinuity}, then all functions from class 
 $p\in{\cal P}_\Omega$ also satisfy the conditions. 

Let us consider a characteristic of function $p$:

$${\cal H}(p)=\sup\limits_{|z|<1}\frac{(1-|z|^2)|p'(z)|}{\Re e p(z)},$$ 
and show that ${\cal H}(p(\varphi(z)))\leqslant {\cal H}(p(z))$ far all $\varphi:\mathbb D\to\mathbb D$.

\begin{equation*}
\begin{split}
	{\cal H}(p(\varphi(z))) = \sup\limits_{|z|<1}
	\frac{(1-|z|^2)| \varphi'(z) p'(\varphi(z))|}
			{\Re e p(\varphi(z))} = \\
	\frac{(1-|z|^2)(1-|\varphi(z)|^2) |\varphi'(z)| |p'(\varphi(z))|}
	{(1-|\varphi(z)|^2)\Re e p(\varphi(z))}.
\end{split}
\end{equation*}
Tacking into account that 

$$
\frac{(1-|z|^2)|\varphi'(z)|}
	{(1-|\varphi(z)|^2)}
\leqslant 1,
$$
we obtain the nessary result.

We have
$$
\frac{\Re z p'(z)}{|z|\Re e p(z)}
\leqslant 
\Big|\frac{p'(z)}{\Re e p(z)}\Big|\leqslant 
\frac{{\cal H}(p)}{(1-r)(1+r)}.
$$
If $\frac 1{1+r}\leqslant 1+\varepsilon$,
then

$$\frac{\Re e z p'(z)}{|z|\Re e p(z)}
	\leqslant 
\frac{{\cal H}(p)(1+\varepsilon)}
{2(1-r)},$$ for $|z|>r$. Taking into account that $\Big|\frac{p'(z)}{\Re e p(z)}\Big|$ is bound in disc $|z|\leqslant r$, 
for all $\varepsilon>0$ we have that
$$
\frac{\Re e z p'(z)}{|z|\Re e p(z)}
	\leqslant 
\frac{{\cal H}(p)(1+\varepsilon)}
{2(1-r)}
	+
O\Big(
\frac{1}
{(1-r)^\alpha}
\Big).$$

Thus to find some result for class  ${\cal P}_\Omega$,  it is sufficient to study only function 
$p_\Omega$.

\begin{theorem}
\label{positivHalfPlan}
If  $\Re e p(z,t) \geqslant k >0$,
then function $p(z,t)$ satisfies conditions of theorem~\ref{inversContinuity}.
\end{theorem}
In this case $\Omega=\{z:\Re e z>k\}$ and function $p_\Omega=(1-k)\frac{1+z}{1-z}+k$.
We only need to show that equation~(\ref{est2}) is true.
Assuming  $z=re^{i\phi}$, thus we have

$$
\frac
	{|p_\Omega|}
	{\Re e p_\Omega}
=\frac
{\Big|  \frac{1+(1-2 k) r e^{i\phi}} {1-r e^{i\phi}}\Big|}
{\Re e\Big( \frac{1+(1-2 k) r e^{i\phi}} {1-r e^{i\phi}} \Big) },0<r<1,\phi\in\mathbb R.$$
From this

$$
\frac
	{|p_\Omega|}
	{\Re e p_\Omega}=
\frac
{ 
  \sqrt{1 + 2 (1 -2 k) r \cos\phi + (1 - 2 k)^2 r^2}
  \sqrt{1 - 2 r \cos\phi+r^2}	
}
{
	1 - 2 k r \cos\phi - (1-2 k) r^2 		 
}. $$

Calculation of derivative on $\psi$ gives us that maximum of this function can be only in points  
$\cos\phi=\pm1,\cos\phi=\frac{2kr}{1-r^2+2kr^2}$. Thus
$$
\frac
	{|p_\Omega|}
	{\Re e p_\Omega}
=
\frac
	{|p_\Omega|}
	{\Re e p_\Omega}
\leqslant 
\max\Big\{
	\sqrt{
		\frac
		{(1+(2k-1)r^2)^2}
		{(1-r^2)(1-(1-2k)r^2)}
	},
	1
\Big\}.$$

Taking into account that $1>k>0$ we have that there is a constant $C_1>0$, such that
$$
\frac
	{|p_\Omega|}
	{\Re e p_\Omega}
\leqslant
C\frac
	1
	{(1-r)^{1/2}
}.
$$

\begin{theorem}
\label{strip}
If $0<a<\Re e p(z,t)<b$  then the solution $w(z,t)$ of equation~(\ref{est2}) is Hoelder function in the
closed unit disk and the boundary of domain the $w(\mathbb D, t)$ is a Jordan curve.
\end{theorem}

In this case $\Omega=\{z:a<\Im m z < b\},0<a<1<b<\infty$ and 
$p_\Omega(w)=\frac{b-a}2\frac{i}{\pi}\log\frac{1+z}{1-z}+\frac{a+b}2.$

$$
{\cal H}(p_\Omega) 
	= 
\sup\limits_{|z|<1}
	\frac{
		(1-|z|^2)
		\frac{b-a}
		{2\pi}
		\frac1{|1-z^2|}
	}
	{
		\frac{b-a}
		{2\pi}
		\arg\frac{1+z}{1-z}
			+
		\frac{a+b}2		
	}.
$$
Assuming $z=re^{i\pi}$ we have
$$
{\cal H}(p_\Omega) 
	= 
\frac{
	(1-r^2)	
}
{	
	|1-r^2e^{i2\psi}|(\arg\frac{1 + r e^{i\psi}}{1 - r e^{i\psi}}
		+
	\delta)		
},
$$
where $\delta = \pi\frac{b+a}{b-a}.$ Thus
$$
{\cal H}(p_\Omega) 
	= 
\frac{
	(1-r^2)	
}
{	
	|1-r^2e^{i2\psi}|(\arg\frac{1 + r e^{i\psi}}{1 - r e^{i\psi}}
		+
	\delta)		
}.
$$
Taking into account that 
$$
\arg
	\frac{1+re^{i\psi}}
	{1-re^{i\psi}}
		=
	\arctan
	\frac{2r\sin\psi}
	{1-r^2}
$$
and
$$
|1 - r^2 e^{i2\psi}| 
	=
\sqrt{
	1 - 2 r^2 \cos 2\psi + r^4
} 
	=
\sqrt{
	(1 - r^2)^2 + 4 r^2 \sin^2 \psi
},$$
we obtain

$$
{\cal H}(p_\Omega) 
	= \frac{
	2
}
{
	\sqrt{1+x^2}
	(\arctan x + \delta)		
},
$$
where $x=\frac{2r\sin\psi}{1-r^2}$.
Using $\sqrt{1+x^2}\arctan x \geqslant 1$ we have

$${\cal H}(p_\Omega)<\frac{1}{1+\delta}.$$

Therefore function  $w(z,t)$  can be extended to a Hoelder function in the
closed unit disk with exponent
$1-\frac{1+\varepsilon}{1+\pi\frac{b+a}{b-a}}$ for each $\varepsilon>0$. Thus boundary $L$ of domain 
$w(\mathbb D,t)$ is curve. From previous theorem we have that inverse function of $w(z,t)$ is continuous on closer of 
domain  $w(\mathbb D,t)$. So each point from 
$L$ has only one preimage. This means that   $L$ is Jordan curve.

\begin{theorem}
\label{sector}
If \begin{equation}
\label{sectorEquation}
\frac{\Im m p(z,t)}{\Re e p(z,t)}<C,
\end{equation}   then the solution $w(z,t)$ of equation~(\ref{est2}) is a Hoelder function in the
closed unit disk and the boundary of domain the $w(\mathbb D, t)$ is a Jordan curve.
\end{theorem}
In this case $\Omega$ is a sector, symmetricalс with respect to real axis, with angle $2\arctan C$.
Thus $p_\Omega=\Big(\frac{1+z}{1-z}\Big)^\alpha,$
where $\alpha=\frac{2}{\pi}\arctan C$.

$$
{\cal H}(p_\Omega)=\frac{\alpha
	\Big|\frac{1+z}{1-z}\Big|^{\alpha-1}\frac{2}{|1-z|^2}
}
{\Re e
	\Big(\frac{1+z}{1-z}\Big)^\alpha
}.
$$

Assuming $z=re^{i\phi}$, we have

$$
{\cal H}(p_\Omega)
=\frac{
	4\alpha\frac{
		|1+re^{i\phi}|^{\alpha-1}
	}
	{
		|1-re^{i\phi}|^{\alpha+1}
	}
}
{\Re e
	\Big(
		\alpha\frac{1+re^{i\phi}}
		     {1-re^{i\phi}}
	\Big)^\alpha +
	\Big(
		\frac{1+re^{-i\phi}}
		{1-re^{-i\phi}}
	\Big)^\alpha
}.
$$
Thus

$$
{\cal H}(p_\Omega)
=\frac{
	4\alpha
	(1+2r\cos\phi+r^2)^{\frac{\alpha-1}2}
	(1-2r\cos\phi+r^2)^{\frac{\alpha-1}2}	
}
{
	\Re e(1+i2r\sin\phi-r^2)^\alpha + (1-i2r\sin\phi-r^2)^\alpha 
}.
$$
Note that

$$\Re (1\pm i2r\sin\phi-r^2)^\alpha\geqslant (1-r^2)^{\alpha},$$
and equality achieves only for $\sin\phi=0$. 
Note that

$$4(1+2r\cos\phi+r^2)^{\frac{\alpha-1}2}(1-2r\cos\phi+r^2)^{\frac{\alpha-1}2}\leqslant 
(1-r^2)^{\frac{\alpha-1}2}$$
and equality achieves only for  $\cos\phi=\pm1$. Thus ${\cal H}(p_\Omega)=2\alpha.$ So,  
function $w(z,t)$ is Hoelder continuous in closed disc with exponent $1-\alpha-\varepsilon$.

From inequality~(\ref{sectorEquation}) follows that $|p_\Omega(z)|\leqslant \sqrt{1+C^2}\Re e p_\Omega(z)$ or conditions 
~(\ref{est1}) in theorem~\ref{inversContinuity}
Taking into account previous argumentation we have
$$\Re e p_\Omega(z) \geqslant 
\frac
{
	(1-r^2)^\alpha
}
{
	(1+r)^{2\alpha}
}=
\frac
{
	(1-r)^\alpha
}
{
	(1+r)^\alpha
},
$$
this means inequality~(\ref{est2}). 
Repeat argumentation from proof of theorem~\ref{strip} we have the necessary result.

\section{Application}

Often in many phenomenas, that are described by Loewner-Kufarev equation,  function
$p(z,t)$  is defined by integral representation.  It is well known, what any function 
$p(z), \Re e p(z) > 0$ can be represented as 

$$
\int\limits_{0}^{2\pi}
\frac{e^{i \theta}+z}
{e^{i \theta}-z}
d\mu(\theta),
$$ 
where  $\mu(\theta)$ is a non-decreasing function.

Let  $\xi(\theta)=\mu'(\theta)$. Thus if $$\xi(\theta)=\frac{1}{|f'(e^{i\theta})|^2},$$ then equation~(\ref{partialLoewnerKufarev}) 
will describe the Hele-Show flow. L.~Carleson and N.~Makarov~\cite{Carleson} used  

$$\xi(\theta)=\inf\{\epsilon:dist(f(e^{i\theta}(1-\epsilon),t),\partial f(\mathbb D,t))=\delta\},\delta > 0$$ 
for the description of processes similar to DLA. In this part 
we connect behaviour of function $\xi(\theta)$ with boundary behaviour of solution of Loewner-Kufarev equation.

First let us assume that $0<a<\xi(\theta)<b$. Then taking into account that, 
$p(0,t) > 0$, we can make a change of time from  $t$ on $\tau$ in equation~\ref{partialLoewnerKufarev} so, that
$p^*(0,\tau)=1$ and $p^*(z,t)=\frac{dt}{d\tau}p(z,t)$. We obtain that
$\frac{a}{b}<p^*(z,\tau)<\frac{b}{a}.$ Thus function $p^*(z,\tau)$ satisfies conditions of theorem~\label{strip}
and we have that solution of Loewner-Kufarev equation maps unit disc onto domain with Hoelder boundary that is Jordan 
curve.

Let us add condition on $\xi(\theta)$, that it is Hoelder continuous this exponent $0<k<1$. Taking into account that (see for 
e.g.~\cite{Mysxelishvili} p.69-79)  

$$|p^*_z(z,t)|=O\Big(\frac{1}{(1-r)^{1-k}}\Big).$$ Note that $\Re p(z,t) > \frac{a}{b}$ 
we have, that the function satisfies conditions of item 2 in theorem~\ref{HolderTheorem}. So we have that the solution $w(z,t)$ 
maps unit disc onto domain with rectifiable and quasiconformal boundary.

\it{Alexander Kuznetsov:}\quad\quad\quad\quad\quad\quad\quad\quad\quad\quad{Alexander.A.Kuznetsov@gmail.com}\newline
\newline
{\it Saratov State University, Department of Mathematics and Mechanics, Astrakhanskaya Str. 83, 410012 Saratov, Russia}
\end{document}